\documentclass[12pt]{amsart}

\usepackage{latexsym}
\usepackage{amsfonts}
\usepackage{amssymb}
\usepackage{fullpage}
\usepackage[OT4]{fontenc}
\usepackage{amscd}
\usepackage{textcomp}

\newtheoremstyle{s2}{9pt}{9pt}{\rm}{\parindent}{\bf}{.}{0.5em}{}
\theoremstyle{s2} 

\newtheorem{definition}{Definition}

\newtheoremstyle{s1}{9pt}{9pt}{\it}{\parindent}{\bf}{.}{0.5em}{}
\theoremstyle{s1}
\newtheorem{lemma}[definition]{Lemma}
\newtheorem{theorem}[definition]{Theorem}
\newtheorem{corollary}[definition]{Corollary}
\newtheorem{proposition}[definition]{Proposition}
\newtheorem{remark}[definition]{Remark}

\newtheorem*{conjecture*}{Conjecture}

\def\P2{{\mathbb{P}^2}}
\def\P1{{\mathbb{P}^1}}
\def\Ce{\mathbb{C}}

\def\Zet{\mathbb{Z}}
\def\En{\mathbb{N}}

\def\Er{\mathbb{R}}

\DeclareMathOperator{\mult}{mult}
\DeclareMathOperator{\lcm}{lcm}
\DeclareMathOperator{\Num}{Num}

\DeclareMathOperator{\gon}{gon}

\numberwithin{definition}{section}

\title[A note on Seshadri constants]
{A note on Seshadri constants of line bundles on hyperelliptic surfaces} \makeatletter

\@addtoreset{equation}{section}
\makeatother
\author{{\L}ucja Farnik}
\address{Jagiellonian University, Institute of Mathematics, {\L}ojasiewicza~6, 30-348 Krak\'{o}w, Poland} \email{lucja.farnik@uj.edu.pl}

\keywords{hyperelliptic surfaces, Seshadri constants, Xu-type lemma.}

\subjclass[2000]{MSC Classification: 14F17, 14E25}

\date{\today}

\begin{document}

\bibliographystyle{alpha}

\begin{abstract}
We study Seshadri constants  of ample line bundles on hyperelliptic surfaces. We obtain new lower bounds and compute the exact values of Seshadri constants in some cases. Our approach uses results of F. Serrano (1990), B. Harboune and J. Ro\'{e} (2008), F.~Bastianelli (2009), A.L. Knutsen, W. Syzdek and T. Szemberg (2009).
\end{abstract}

\maketitle
\thispagestyle{empty}

\section{Introduction}
Seshadri constants measure how positive a line bundle is. They were introduced in 1992 by  J.P.~Demailly in \cite{De1992} as an attempt to tackle the famous Fujita conjecture. The conjecture has not been proven but Seshadri constants soon became an object of study on their own.

Giving exact values or just estimating Seshadri constants is very hard, even in case of line bundles on algebraic surfaces, see e.g. \cite{Ba1999}. There exists an upper bound for Seshadri constant of a line bundle at points $x_1$, $\ldots$, $x_r$ on a smooth projective $n$-dimensional
variety~$X$, namely $\varepsilon(L,x_1, \ldots, x_r)\leq \sqrt[n]{\frac{L^n}{r}}$. 
Therefore it is interesting to look for lower bounds.

There are several results concerning Seshadri constants on surfaces with Kodaira dimension zero. Let us recall some of them. In  appendix to \cite{Ba1998}, Th. Bauer and T. Szemberg give upper bound for the global Seshadri constant of an ample line bundle on an abelian surface and as a corollary obtain that the Seshadri constant of such a line bundle is always rational. In \cite{Ba1997} Th. Bauer computes Seshadri constants on all $K3$ surfaces of degree 4. This result is extended by 
 C. Galati and A.L. Knutsen in \cite{GaK2013} who compute Seshadri constants  on $K3$ surfaces of degrees 6 and 8. Earlier in \cite{K2008}  A.L. Knutsen estimates Seshadri constants on $K3$ surfaces with Picard number 1.
T. Szemberg in \cite{Sz2001} proves that the global Seshadri constants on  Enriques surfaces are always rational and also provides the lower bound for Seshadri constant at an arbitrary point.
Up to our knowledge, Seshadri constants have not been studied on hyperelliptic surfaces before.

We estimate Seshadri constants on hyperelliptic surfaces, in some cases we compute their exact values. 
The paper is organised in the following way --- in Theorem   \ref{StSesh1nahiper} we compute the global Seshadri constant of a line bundle of type (1,1) on a hyperelliptic surface of an arbitrary type, next in Proposition  \ref{epislon>1} we point out a hyperelliptic surface type and a point at which a Seshadri constant of  a line bundle of type (1,1) is strictly greater than   $1$. 
In Theorem  \ref{stSesh_wiazki_szer_na_pow_1} we compute the global Seshadri constant of a arbitrary ample line bundle on hyperelliptic surface of type $1$, and in Theorem ~\ref{stSesh_wiazki_szer_na_pow_2-7} we provide a lower bound for this constant on  hyperelliptic surfaces of types $2$-$7$.
Finally, in Theorem \ref{HRdlahiper} we  estimate from below the multi-point Seshadri constant of an ample line at $r$ very general points on hyperelliptic surfaces.

\section{Notation and auxiliary results}

Let us set up the notation and basic definitions. We work over the field of complex numbers~$\Ce$. We consider only smooth reduced and irreducible projective varieties. By $D_1\equiv D_2$ we denote the numerical equivalence of divisors $D_1$ and $D_2$. By a curve we understand an irreducible subvariety of dimension 1.
In the notation we follow \cite{Laz2004}.

Let $X$ be a smooth projective
variety and $L$ a nef line bundle on~$X$. 
We recall the definition of a Seshadri constant. 

\begin{definition} 
(1)
The Seshadri constant of $L$ at a given point $x\in X$ is the real number
$$\varepsilon(L,x)=\inf \left\{\frac{LC}{\mult_xC}: \  C\ni x\right\},$$
where the infimum is taken over all irreducible curves 
$C\subset X$ passing through $x$.

(2) The global Seshadri constant of $L$ is defined to be 
$$\varepsilon(L)=\inf_{x\in X} \varepsilon(L,x).$$
\end{definition}

Let  $x_1$, $\ldots$, $x_r$ be pairwise distinct points. The notion of a Seshadri constant of a line bundle at a point may be generalised to $r$ points in the following way:
\begin{definition}
 The {multi-point Seshadri constant} of $L$ at $x_1$, $\ldots$, $x_r$
is the real number
$$\varepsilon(L,x_1, \ldots, x_r)=\inf \left\{\frac{LC}{\sum_{i=1}^r\mult_{x_i}C}: \ \{x_1,\ldots,x_r\}\cap C\neq \emptyset\right\},$$
where the infimum is taken over all irreducible curves $C\subset X$ passing through at least one of the points 
 $x_1$, $\ldots$, $x_r$.
\end{definition}

For a fixed line bundle $L$ the function  $(x_1, \ldots, x_r)\mapsto \varepsilon(L,x_1, \ldots, x_r)$ is constant for points in very general position, moreover its value for points not lying in very general position does not exceed the value for points in very general position --- see \cite{Laz2004}, Example 5.1.11.
We denote the Seshadri constant of $L$ at $r$ points in very general position by  $\varepsilon(L,r)$.

 Let $\alpha_0(L, m_1,\ldots, m_r)$ denote the least degree $LC$ of an irreducible curve $C$ passing through $r$ points in general position with multiplicities $m_1,\ldots, m_r$. Let $m^{[l]}=\underbrace{(m, \ldots, m)}_{l \text{ times}}$. Then the following theorem holds:

\begin{theorem}[Harbourne, Ro\'{e}, \cite{HR2008}, Theorem 1.2.1]\label{HR2008}
Let $L$ be a big and nef line bundle on a smooth projective surface. Let $r\in\En$, $r\geq 2$, let $\mu\in\Er$, $\mu \geq 1$. If

(1) for every $m\in\En$ such that $1\leq m<\mu$ we have
$$\alpha_0(L,m^{[r]})\geq m\sqrt{L^2\left(r-\frac{1}{\mu}\right)}$$

\hskip 20pt  and

(2) for every $m\in\En$ such that $1\leq m<\frac{\mu}{r-1}$ and if for every $k\in\Zet$ such that $k^2<\frac{r}{r-1}\min\{m, m+k\}$ we have
$$\alpha_0(L,m^{[r-1]},m+k)\geq \frac{mr+k}{r}\sqrt{L^2\left(r-\frac{1}{\mu}\right)},$$

\hskip 20pt then 
$$\varepsilon(L, r)\geq\sqrt{\frac{L^2}{r}}\sqrt{1-\frac{1}{r\mu}}.$$
\end{theorem}

For more background on Seshadri constants we refer to \cite{PSC2009}.

\vskip 10pt
Now let us recall the definition of a hyperelliptic surface.

\begin{definition}
A hyperelliptic surface $S$ (sometimes called bielliptic) 
is a surface with Kodaira  dimension equal to $0$ and  irregularity  $q(S)=1$.
\end{definition}

Alternatively (\cite{Bea1996}, Definition VI.19), a surface $S$ is hyperelliptic if  $S\cong (A\times B)/G$, where $A$ and $B$ are elliptic curves, and $G$ is an abelian group acting on A by translation and acting on B, such that  $A/G$ is an elliptic curve and $B/G\cong \mathbb{P}^1$; $G$ acts on $A\times B$ coordinatewise.
Hence we have the following situation:
$$
\begin{CD}
S\cong (A\times B)/G @>\Phi>> A/G @.\\
@V\Psi VV @.\\
B/G\cong \P1
\end{CD}
$$
where $\Phi$ and $\Psi$ are natural projections.

Hyperelliptic surfaces were classified at the beginning of 20th century by G. Bagnera and M.~de Franchis in \cite{BF1907}, and independently by F. Enriques i F. Severi in \cite{ES1909-10}. They showed that there are seven non-isomorphic types 
of hyperelliptic surfaces. Those types are characterised by the action of $G$ on $B\cong\Ce/(\Zet\omega\oplus\Zet)$ (for details see eg. \cite{Bea1996}, VI.20). For every hyperelliptic surface we have that the canonical divisor $K_S$ is numerically trivial.

In 1990 F.~Serrano in \cite{Se1990} characterised the group $\Num(S)$
for each of the  surface's type:
\begin{theorem}[Serrano] A basis of the group of classes of numerically equivalent divisors  $\Num(S)$
for each of the  surface's type and the multiplicities of the singular fibres in each case are the following:
$$
\begin{array}{c|l|l|l}
\text{Type of a hyperelliptic surface}&G&m_1,\ldots,m_s&\text{Basis of $\Num(S)$}\\
\hline
1&\Zet_2&2,2,2,2&A/2, B\\
2&\Zet_2\times\Zet_2&2,2,2,2&A/2, B/2\\
3&\Zet_4&2,4,4&A/4, B\\
4&\Zet_4\times\Zet_2&2,4,4&A/4, B/2\\
5&\Zet_3&3,3,3&A/3, B\\
6&\Zet_3\times\Zet_3&3,3,3&A/3, B/3\\
7&\Zet_6&2,3,6&A/6, B
\end{array}
$$
\end{theorem}
Let $\mu=\lcm\{m_1, \ldots, m_s\}$ and let $\gamma=|G|$. Given a hyperelliptic surface, its basis of  $\Num(S)$ consists of divisors $A/\mu$ and $\left(\mu/\gamma\right) B$. We say that $L$ is a line bundle of type $(a,b)$ on a hyperelliptic surface if $L\equiv a\cdot A/\mu+b\cdot(\mu/\gamma) B$.
In  $\Num(S)$  we have that $A^2=0$, $B^2=0$, $AB=\gamma$. 
Note that a divisor  $b\cdot\left(\mu/\gamma\right) B\equiv (0,b)$, $b\in\Zet$, is effective if and only if  $b\cdot\left(\mu/\gamma\right)\in \En$ (see \cite{Ap1998}, Proposition~5.2).

The following proposition holds:
\begin{proposition}[see \cite{Se1990}, Lemma 1.3]\label{Ser2}
Let $D$ be a divisor of type $(a,b)$ on a hyperelliptic surface $S$. Then
\begin{enumerate}
\item  $\chi(D)=ab$;
\item $D$  is ample if and only if $a>0$ and  $b>0$;
\item If $D$  is ample then $h^0(D)=\chi(D)=ab$.
\end{enumerate}
\end{proposition}

\vskip 5pt
Now we recall a bound for the self-intersection of a curve.
Adjunction formula, applied to the normalisation of a curve $C$, implies the following formula:
\begin{remark}[Genus formula, \cite{GH1978}, Lemma, p. 505] Let $C$ be a curve on a surface $S$, passing through $x_1$, $\ldots$, $x_r$ with multiplicities respectively $m_1$, $\ldots$, $m_r$. Let $g(C)$ denote the genus of the normalisation of $C$.
Then
$$g(C)\leq \frac{C^2+C.K_{S}}{2}+1-\sum_{i=1}^r\frac{m_i(m_i-1)}{2}.$$
\end{remark}

Note that a curve $C$ on a hyperelliptic surface has  genus at least 
 $1$. Otherwise the normalisation of $C$, of genus zero, would be a covering (via  $\Phi$) of an elliptic curve  $A/G$. This contradicts the Riemann-Hurwitz formula.

\vskip 5pt

For families of curves we have Xu-type lemma.
The original version of this lemma was proved by  G.~Xu in \cite{Xu1995}. We will use the generalisation of the Xu Lemma obtained by  A.L.~Knutsen, W.~Syzdek, T.~Szemberg in  \cite{KSSz2009}, and independently by F. Bastianelli in \cite{Bas2009}. 
Let $\gon({C})$ denote the gonality  of a smooth  curve  $C$, i.e. the minimal degree of a covering ${C}\to \P1$.

\begin{lemma}[Bastianelli, \cite{Bas2009}, Lemma 2.2; Knutsen-Syzdek-Szemberg, \cite{KSSz2009}, Theorem A]\label{Xu-KSSz}
Let $S$ be a smooth projective surface. Let $U$ be a smooth variety. Consider a nontrivial family $\{(C_u,x_u)\}_{u\in U}$ where $x_u$ is a very general point of $S$ and  $C_u$ is a curve satisfying  the condition $\mult_{x_u}C_u \geq m$ for every  $u\in U$ and for some integer $m\geq 2$. Then for a general curve $C$ of this family 
$$C^2 \geq m(m-1)+\gon(\widetilde{C}).$$
\end{lemma}

Applying the Xu-type lemma to a family $\mathcal{C}$ of curves passing through $x_1$, $\ldots$, $x_r$ with multiplicities respectively $m_1$, $\ldots$, $m_r$, where $m_1\geq 2$, on a blow-up at $x_2$, $\ldots$, $x_{r}$,
 we have the following multi-point version of the Xu-type lemma
\begin{lemma}\label{Xu-KSSz2}
For a general curve $C$ of the family ${\mathcal{C}}$ as above we have that
$$C^2 \geq \left(\sum_{i=1}^r m_i^2\right)-m_1+\gon(\widetilde{C}).$$
\end{lemma}

Every hyperelliptic surface $S$ is nonrational, hence for every curve $C\subset S$ we have \mbox{$\gon(\widetilde{C})\geq 2$} (see \cite{KSSz2009}, remarks following Theorem A).

\section{Main results}

\subsection{Seshadri constants of ample line bundles on hyperelliptic surfaces}

We start with computing the global Seshadri constant in the simplest case of an ample line bundle on a hyperelliptic surface, i.e. for a line bundle of type $(1,1)$.
\begin{theorem}\label{StSesh1nahiper}
Let $S$ be a hyperelliptic surface. Let $L$ be a line bundle of type $(1, 1)$ on $S$. Then $$\varepsilon(L)= 1.$$

\begin{proof}
Let $C\equiv(\alpha,\beta)$ denote an irreducible curve passing through a given point $x\in S$ with multiplicity $m$, $m\geq 1$. We estimate the value of  $\frac{LC}m$ from below.

Depending on the position of a point $x$ and on hyperelliptic surface's type, we have the following possibilities for $C$ to be an irreducible curve:

(1) $C\equiv B\equiv(0,k)$ and $x$ is an arbitrary point, where $k=1$ for  a hyperelliptic surface of an odd type; $k=2$ for  a hyperelliptic surface of type 2 and 4; $k=3$ for  a hyperelliptic surface of type 6. Then
$$\frac{LC}m=\frac k 1\geq 1.$$

(2) $C\equiv nA/\mu\equiv(n,0)$ and a point $x$ lies on a fibre $nA/\mu$,  where $n\in\{1,2\}$ for  a hyperelliptic surface of type $1$ and $2$;  $n\in\{1,2,4\}$ for type $3$ and $4$; $n\in\{1,3\}$ for type $5$ and~$6$; $n\in\{1,2,3,6\}$  for type $7$. Then $$\frac{LC}m=\frac n 1\geq 1.$$

(3) $C\equiv(\alpha, \beta)$, where $\alpha> 0$ and $\beta> 0$, and $x$ is an arbitrary point. Then by B{\'e}zout's Theorem, intersecting  $C$ with a fibre $B$ and with an appropriate, depending on the position of the point  $x$, fibre $nA/\mu$, we get:
$$\frac{LC}m=\frac{\alpha+\beta}m\geq\left\{\begin{array}{ll}
1,&\text{in case of a hyperelliptic surface of type } 1, 3,5,7;\\
\frac 1 2+\frac 1 2,&\text{in case of a hyperelliptic surface of type } 2;\\
\frac 1 2+\frac 1 4,&\text{in case of a hyperelliptic surface of type } 4;\\
\frac 1 3+\frac 1 3,&\text{in case of a hyperelliptic surface of type } 6.
\end{array}\right.$$

\noindent
Therefore $\frac{LC}m\geq 1$ for a hyperelliptic  surface of type $1$, $2$, $3$, $5$ and $7$. 

Now let $S$ be a surface of type $4$ or $6$. We consider two cases. If $m=1$, then $\frac{LC}{m}=\frac{\alpha+\beta}{m}\geq \frac{2}{1}>1$. If $m\geq 2$, then by genus formula $C^2\geq  m^2-m$, by Hodge Index Theorem 
$(LC)^2\geq L^2C^2= 2C^2\geq 2(m^2-m)$, and therefore
$\frac{LC}{m}\geq \sqrt{\frac{2(m^2-m)}{m^2}}=\sqrt{2-\frac{2}{m}}\geq 1$. 

Hence independently of the type of the hyperelliptic surface we have $\varepsilon(L,x)\geq 1$.
 Moreover, for every hyperelliptic surface's type, $\varepsilon(L,x)= 1$ for a point $x$ on a fibre $A/\mu$.
 Therefore $\varepsilon(L)= 1$.
\end{proof}
\end{theorem}

By the proof of Theorem \ref{StSesh1nahiper} we immediately obtain a corollary:
\begin{corollary}
Let $S$ be a hyperelliptic surface of an odd type. Let $L$ be a line bundle of type $(1, 1)$ on $S$. Then the Seshadri constant of $L$ at any $x\in S$ is computed by a fibre $B$, hence $$\varepsilon(L,x)= 1 \text{ for any } x\in S.$$
\end{corollary}

On the other hand, it is not true that on every hyperelliptic surface the equality  $\varepsilon(L,x)= 1$ holds  for every  $x\in S$.

\begin{proposition}\label{epislon>1}
There exists a hyperelliptic surface $S$ such that for a line bundle $L$ of type $(1,1)$
$$\varepsilon(L,1)>1.$$
\begin{proof}
Let  $S$ be a hyperelliptic surface of type  $2$, and let $L$ be a line bundle of type $(1,1)$ on~$S$. Let $x$ be a very general point on $S$. 
We will prove that $\varepsilon(L,x)\geq \frac 4 3$.

Let $C\equiv(\alpha,\beta)$ be an irreducible curve passing through a given point $x\in S$ with multiplicity $m$, $m\geq 1$.

Let $m=1$. Assume that  $\frac{LC}m<\frac 4 3$. Then $LC<\frac 4 3$, hence $\alpha+\beta<\frac 4 3$ and thus, as $\alpha$ and $\beta$ are nonnegative integers,  $\alpha+\beta\leq 1$. Therefore either  $(\alpha,\beta)\equiv(1,0)\equiv  A/2$, or $(\alpha,\beta)\equiv(0,1)\equiv B/2$. A divisor $A/2$ does not pass through $x$, a divisor  $B/2$ is not effective, a contradiction.

Now let $m\geq 2$. We have to prove that $\frac{LC}m\geq\frac 4 3$. 
Both sides are nonnegative, hence equivalently $(LC)^2\geq\frac{16}{9} m^2$. By Hodge Index Theorem it is enough to show that
$L^2C^2\geq\frac{16}{9} m^2$. By Xu-type lemma (Lemma \ref{Xu-KSSz}) we have 
that $C^2\geq m^2-m+2$. Hence it is enough to prove that
 $2m^2-2m+4\geq \frac{16}{9} m^2$. Equivalently $(m-3)(m-6)\geq 0$.
The inequality is satisfied for $m\neq 4,5$. We consider these two cases separately.

Let $m=4$. Suppose that  $\frac{LC}4<\frac 4 3$. Hence ${LC}<\frac {16} {3}$, so $\alpha+\beta\leq 5$. On the other hand, by Xu-type lemma $2\alpha\beta=C^2\geq m^2-m+2=14$, a contradiction.

For $m=5$, if $\frac{LC}5<\frac 4 3$ then $\alpha+\beta\leq 6$. By Xu-type lemma $\alpha\beta\geq 11$,  a contradiction. This completes the proof.
\end{proof}
\end{proposition}

Using the same method as presented in Theorem \ref{epislon>1} one can show that for a very general point $x$ on a hyperelliptic surface of type 2 and  for $L$ of type $(1,1)$, the Seshadri constant of $L$ at $x$ is greater than a constant slightly bigger than $\frac 4 3$. The proof splits in a large number of cases and therefore we decide not to present it here. However precise study of this example might support the idea that
this Seshadri constant is irrational.

\vskip 5pt
Now we will prove a lower bound for the global Seshadri constant of an arbitrary ample line bundle on hyperelliptic surface of type 1.

\begin{theorem}\label{stSesh_wiazki_szer_na_pow_1}
Let $S$ be a hyperelliptic surface of type 1. Let $L$ be an ample line bundle of type $(a,b)$ on $S$. Then $$\varepsilon(L)=\min \{a,b\}.$$

\begin{proof}
Let $S$ be a hyperelliptic surface of type 1, let $L\equiv (a,b)$. 
Let $C\equiv(\alpha,\beta)$ denote an irreducible curve passing through a given point $x$ with multiplicity $m$, $m\geq 1$. Using B{\'e}zout's Theorem we obtain:
$$\frac{LC}m=\frac{a\beta+b\alpha}m\geq\left\{\begin{array}{ll}
a,&\text{if } C\equiv B \text{ and } x \text{ is an arbitrary point};\\
b,&\text{if } C\equiv A/2 \text{ and } x \text{ lies on the singular fibre } A/2;\\
2b,&\text{if } C\equiv A \text{ and } x \text{ lies on the fibre } A;\\
a+b,&\text{if } C\equiv (\alpha,\beta) \text{ and } x \text{ lies on one of the singular fibres } A/2;\\
\frac a 2+b,&\text{if } C\equiv (\alpha,\beta) \text{ and } x \text{ lies on one of the  fibres } A.
\end{array}\right.$$

Hence on a hyperelliptic surface of type 1  
\begin{equation*}
\varepsilon(L)= \min\{a,b\}. \qedhere
\end{equation*}

\end{proof}
\end{theorem}

By the theorem above we see that on a hyperelliptic surface of type 1 the global Seshadri constant of an ample line bundle $L$ is always submaximal, ie. smaller than $\sqrt{L^2}$.

Note that the method used in Theorem  \ref{stSesh_wiazki_szer_na_pow_1} does not work on hyperelliptic surfaces of other types. For hyperelliptic surfaces of type 1 the lower bound of  $\frac{LC}m$, where a curve $C$ is not a  fibre, is always greater than the value of $\frac{LC}m$ for some fibre~$C$. It is also easy to show for which fibre and for  which point position the global Seshadri constant is actually reached.  This is not the case for hyperelliptic surfaces of types 2-7.

\vskip 5pt
For hyperelliptic surfaces of types 2-7, we have the following lower bound for the global Seshadri constant
\begin{theorem}\label{stSesh_wiazki_szer_na_pow_2-7}
Let $S$ be a hyperelliptic surface of type greater than 1. Let $L$ be an ample line bundle of type $(a,b)$ on $S$. Then $$\varepsilon(L)\geq \min\{a,b\}.$$
\begin{proof}
We have that $L\equiv (a,b)\equiv \min\{a,b\}\cdot M+N$, where $M\equiv (1,1)$ and $N$ is nef. By definition of a Seshadri constant, for every $x\in S$
$$\varepsilon(L,x)\geq \min\{a,b\}\cdot\varepsilon(M,x)+\varepsilon(N,x)\geq \min\{a,b\}\cdot\varepsilon(M,x).$$ Hence by Theorem \ref{StSesh1nahiper}
\begin{equation*}
\varepsilon(L)\geq \min\{a,b\}\cdot\varepsilon(M)=\min\{a,b\}.\qedhere
\end{equation*}
\end{proof}
\end{theorem}
\vskip 5pt

\subsection{Multi-point Seshadri constants of ample line bundles on non-rational surfaces}
In this section we present a lower bound for  Seshadri constant at $r$ points in very general position on  hyperelliptic surfaces.

The lower bound for multi-point Seshadri constants obtained in Theorem \ref{HRdlahiper} is not far from the upper bound. As mentioned before, it is well-known (see e.g. \cite{PSC2009}, Proposition~2.1.1) that for smooth projective surfaces $$\varepsilon(L, r)\leq\sqrt{\frac{L^2}{r}}.$$

Biran-Nagata-Szemberg conjecture says that for any algebraic surface there exists $r_0>0$ such that for every $r>r_0$ in fact there is an equality $\varepsilon(L, r)=\sqrt{\frac{L^2}{r}}.$

\begin{theorem}\label{HRdlahiper}
Let $S$ be a hyperelliptic surfaces. Let $L$ be an ample line bundle on $S$. Then
\begin{equation*}\label{HSesh-r}
\varepsilon(L, r)\geq\sqrt{\frac{L^2}{r}}\sqrt{1-\frac{1}{8r}}, \quad r\geq 2.
\end{equation*}

\begin{proof}
The claim follows immediately from Harbourne-Ro\'{e} theorem (Theorem \ref{HR2008}) with $\mu=8$. The point is to check that the assumptions of this theorem are satisfied with that particular constant. Turning into details we need to check the following two conditions: 

(1) for every integer  $1\leq m<8$
$$\alpha_0(L,m^{[r]})\geq m\sqrt{L^2\left(r-\frac{1}{8}\right)};$$

(2)  for every integer $1\leq m<\frac{8}{r-1}$ and  for every integer $k$  with $k^2<\frac{r}{r-1}\min\{m, m+k\}$
$$\alpha_0(L,m^{[r-1]},m+k)\geq \frac{mr+k}{r}\sqrt{L^2\left(r-\frac{1}{8}\right)}.$$

\vskip 5 pt

Ad. (1). For $m=1$, $2$, $\ldots$, $7$ we ask whether the inequality
$$\alpha_0(L,m^{[r]})  \geq  m\sqrt{L^2\left(r-\frac{1}{8}\right)}$$ is satisfied.

Let $C$ be an irreducible curve computing $\alpha(L,m^{[r]})$. It suffices to show that
$$LC \geq    m\sqrt{L^2\left(r-\frac{1}{8}\right)}.$$

As $L$ is ample, by Hogde Index Theorem it is enough to prove that
$$L^2C^2 \geq   m^2 {L^2\left(r-\frac{1}{8}\right)}.$$

We split the proof that 
$C^2 \geq   m^2 {\left(r-\frac{1}{8}\right)}$ into
 two cases:  $m=1$ and $m>1$.

For $m=1$, we have  $h^0(C)=\dim|C|+1\geq r+1$. Moreover by Proposition \ref{Ser2} (3), $h^0(C)=\frac{C^2}{2}$. Hence $\frac{C^2}{2}\geq r+1$. Therefore it is enough to show that
$$2r+2\geq     r-\frac{1}{8}.$$
This condition is satisfied for every positive $r$.

\vskip 5pt

Now let $2\leq m\leq 7$. By Xu-type lemma (Lemma \ref{Xu-KSSz2}), 
$C^2\geq r m^2 - m+2$. Hence it is enough to show that
$$r m^2 - m+2\geq    m^2 {\left(r-\frac{1}{8}\right)},$$
which is elementary.

\vskip 5pt

Ad. (2). In the table below we write down all  values of $r$, $m$ and  $k$ satisfying  conditions
  $1\leq m<\frac{8}{r-1}$ and $k^2<\frac{r}{r-1}\min\{m, m+k\}$.
$$\begin{array}{l | c |l}
r&m<\frac{8}{r-1}&\text{possible } k\\
\hline
2&1&1\\
&2&1,-1\\
&3&1,-1,2\\
&4&1,-1,2\\
&5&1,-1,2,-2,3\\
&6&1,-1,2,-2,3\\
&7&1,-1,2,-2,3\\
\hline
3&1&1\\
&2&1,-1\\
&3&1,-1,2\\
\hline
4&1&1\\
&2&1,-1\\
\hline
5&1&1\\
\hline
6&1&1\\
\hline
7&1&1\\
\hline
8&1&1\\
\end{array}$$
We have omitted the case $k=0$ in each row, as for $k=0$ we get the inequality already proved in~(1).

Using Hodge Index Theorem, analogously to (1) the condition to check is reduced to the following inequality
$$C^2\geq    \left(\frac{mr+k}{r}\right)^2\left(r-\frac{1}{8}\right), $$
where $C$ be an irreducible curve computing $\alpha_0(L,m^{[r-1]},m+k)$.

Again we consider two cases: $m=1$ and {$m>1$}.

Let $m=1$. Hence $k=1$. Since Xu-type lemma (Lemma \ref{Xu-KSSz2}) implies that $C^2\geq r+3$, we easily obtain that
$$C^2 \geq    \left(\frac{r+1}{r}\right)^2\left(r-\frac{1}{8}\right).$$

\vskip 5pt

Let $m>1$. Hence $r\in\{2,3,4\}$. By Xu-type lemma $C^2\geq\ (r-1) m^2 + (m+k)^2-m+2$, so it is enough to show that
$$ (r-1) m^2 + (m+k)^2-m+2\geq  \left(\frac{mr+k}{r}\right)^2\cdot\left(r-\frac{1}{8}\right)$$
holds. After reordering the terms we obtain the following inequality
$$8r^2k^2-8r^2m+16r^2+m^2r^2+2mrk-8rk^2+k^2\geq  0.$$
Simple computations confirm that the last inequality is satisfied for all admissible  $m>1$, $r$ and $k$. The proof is completed.
\end{proof}
\end{theorem}

\begin{remark}
Note that the Theorem \ref{HSesh-r} holds  also for abelian surfaces with $\rho=1$.
\end{remark}

\subsection*{Acknowledegments} The author would like to thank Halszka Tutaj-Gasi\'nska for advice and support, and Tomasz Szemberg for many helpful discussions and improving the readability of the paper.

\end{document}